\documentclass[11]{amsart}
\usepackage{amsfonts,amsmath,amssymb,eulervm}
\usepackage[all]{xy}
\usepackage{hyperref}

\begin{document}

\newcommand{\ba}{{\bf a}}
\newcommand{\cra}{{\boldsymbol{\alpha}}}
\newcommand{\bb}{{\bf b}}
\newcommand{\bB}{{\sf B}}
\newcommand{\fd}{{\bf d}}
\newcommand{\bx}{{\bf x}}
\newcommand{\ab}{{\rm ab}}
\newcommand{\bpi}{{\boldsymbol{\pi}}}
\newcommand{\be}{{\bf e}}
\newcommand{\oh}{{\mathfrak o}}
\newcommand{\m}{{\mathfrak m}}
\newcommand{\jnf}{{\rm inf}}
\newcommand{\A}{{\mathbb A}}
\newcommand{\cA}{{\mathcal A}}
\newcommand{\B}{{\mathbb B}}
\newcommand{\cB}{{\mathcal B}}
\newcommand{\C}{{\mathbb C}}
\newcommand{\E}{{\mathcal E}
}
\newcommand{\F}{{\mathbb F}}
\newcommand{\D}{{\mathbb D}}
\newcommand{\G}{{\mathbb G}}
\newcommand{\fH}{{\mathfrak H}}
\newcommand{\cH}{{\mathcal H}}
\newcommand{\fh}{{\mathfrak h}}
\newcommand{\M}{{\mathbb M}}
\newcommand{\N}{{\mathbb N}}
\newcommand{\bO}{{\boldsymbol{\Omega}}}
\newcommand{\bP}{{\mathbb P}}
\newcommand{\R}{{\mathbb R}}
\newcommand{\bS}{{\mathbb S}}
\newcommand{\SU}{{\rm SU}}
\newcommand{\Gal}{{\rm Gal}}
\newcommand{\Q}{{\mathbb Q}}
\newcommand{\bQ}{{\overline{\mathbb Q}}}
\newcommand{\T}{{\mathbb T}}
\newcommand{\Sp}{{\rm Sp \;}}
\newcommand{\Spin}{{\rm Spin}}
\newcommand{\W}{{\sf W}}
\newcommand{\ga}{{\mathfrak a}}
\newcommand{\sK}{{\sf K}}
\newcommand{\Z}{{\mathbb Z}} 
\newcommand{\Be}{{\boldsymbol{\epsilon}}}
\newcommand{\Maps}{{\rm Maps}}
\newcommand{\Mor}{{\rm Mor}}
\newcommand{\Stab}{{\rm Stab}}
\newcommand{\sk}{{\rm sk}}

\newcommand{\weta}{{\vartheta}}
\newcommand{\sslash}{{/\!/}}

\newcommand{\Spec}{{\rm Spec \; }}
\newcommand{\kerr}{{ \{ \rm Kerr \} }}
\newcommand{\tz}{{\vartheta}}

\parindent=0pt
\parskip=6pt

\newcommand{\ie}{\textit{ie}\,}
\newcommand{\eg}{\textit{eg}\,}
\newcommand{\se}{{\sf e}}
\newcommand{\cf}{{\textit{cf}\,}}
\newcommand{\K}{{\rm K}}
\newcommand{\Rep}{{\rm R}}
\newcommand{\ox}{{\rm or}}
\newcommand{\bM}{{\overline{\mathbb M}}}
\newcommand{\barpi}{{\boldsymbol {\varpi}}}
\newcommand{\bE}{{\bf \D/\cpct_0}}
\newcommand{\bH}{{\mathbb H}}
\newcommand{\bl}{{\l{}}}
\newcommand{\wzw}{{\bf wzw}}
\newcommand{\Pic}{{\rm Pic}}
\newcommand{\ctc}{{\sf Ctct}}
\newcommand{\cpct}{{\rm Cpct}}
\newcommand{\lcs}{{\sf lcs}}
\newcommand{\bL}{{\boldsymbol{\Lambda}}}
\newcommand{\brh}{{\boldsymbol{\rho}}}
\newcommand{\Mod}{{\rm Mod}}

\newcommand{\OBlcs}{{\bO^*B \lcs}}
\newcommand{\qt}{{_{\Q[t]}}}
\newcommand{\OS}{{\bO^*\Omega S^2}}
\newcommand{\Fin}{{\rm Fin}}
\newcommand{\SL}{{{\mathfrak s}{\rm l}}}
\newcommand{\Met}{{\M{\rm et}}}
\newcommand{\iso}{{\rm iso}}
\newcommand{\HAut}{{{\mathcal H}{\rm Aut}}}
\newcommand{\halpha}{{\cancel \alpha}}
\newcommand{\cS}{{\mathcal S}}
\newcommand{\pt}{{\rm pt}}
\newcommand{\HL}{{{\sf HL}}}
\newcommand{\Vect}{{{\rm Vect}_\R}}
\newcommand{\Fix}{{\rm Fix}}
\newcommand{\Gl}{{{\mathbb G}{\rm l}}}
\newcommand{\oso}{{\rm Iso}}
\newcommand{\AQ}{{\mathcal{AQ}}}
\newcommand{\cX}{{\mathcal X}}
\newcommand{\cY}{{\mathcal Y}}


\title{Boundary framings for $\lcs$ four-manifolds}

\author[J Morava]{J Morava}

\address{Department of Mathematics, The Johns Hopkins University,
Baltimore, Maryland} 

\email{jmorava1@jhu.edu}

\begin{abstract}{In work with Urs Schreiber [\S 1] we construct a rational homotopy-theoretic model for a classifying space of locally conformally symplectic ($\lcs$) structures on four-manifolds. This leads to the definition of a cobordism category of three-manifolds `anchored' by principal $\Omega^2 S^2$ - bundles (\S 2, generalizing contact structures).  Powerful $\SL_2\R$ - representation-valued Hodge-Lefschetz cohomology functors (\cite{3, 8, 41, 42}, going back to Chern and Weil), taking values in a category $\cB^*$ of bidifferential $\Z$-graded complexes \cite{2}(lemma 1.3), are available for its study. \medskip

\noindent \S3 observes that the $\Q$-rationalization of the Ebin - Omori \cite{9, 36} configuration groupoid $[\Met \sslash \D]$ for Euclidean general relativity can be presented in terms of a category of elements of a 4D manifold $X$ under the action of its group $\D$ of diffeomorphisms. Thus is essentially a notation system for actions of compact connected groups $K$ on $X$ respecting various geometric structures. This revision couples topological eGR to $\lcs$-Hamiltonian matter via the forgetful functor 
\[
[\lcs \times_{\G_m} \Met \sslash \cS] \Rightarrow [\lcs \times_{\G_m} \Met \sslash \D] {\;}.
\]  
}\end{abstract}

\maketitle 

{\sc TOC}

{\bf \S I}  A homotopy-theoretic classifying space for $\lcs$ structures:\\
(including contact and open book structures, see also {\bf a2})

{\bf \S 2} Spin, contact, and $\bl$ine structures:\\
3-manifolds with $\bl$ine structures as objects of a cobordism category with $\lcs$ manifolds as morphisms

{\bf \S 3} Delayed introduction : {\bf Why might this be interesting?}

{\bf appendices}

{\bf 1} contact structures and open books\\
{\bf 2} Hodge-Lefschetz cohomology of $\lcs$ manifolds\\
{\bf 3} details re \S 2.1\\
{\bf 4} Spaces $X \to |\Fin^\otimes_\C|^+$ (with enough points)\\
{\bf 5} details from [HR] re flux\\
{\bf 6} Kerr spacetime

\newpage

{\bf \S I  A homotopy-theoretic classifying space for $\lcs$ structures} 

{\bf 1.1} In the following, $ X \supset \partial X = Y$ will be a connected compact oriented $2n$-dimensional manifold with (for now) connected boundary, together with a locally conformally symplectic ($\lcs$) structure \cite{15, 16, 17, 21} with an atlas of compatible two-forms $\omega \in \bO^2(U), U \subset X$ defining local volumes $\omega^n/n! \neq 0$, together with their unique closed \cite{40} one-forms $\weta \in \bO^1(U)$ such that
\[
d\omega = \weta \wedge \omega \;.
\]
A morphism $(\phi,u) : (X',\omega') \to (X,\omega)$ of $\lcs$ manifolds (\eg charts) is a diffeomorphism $\phi$ with
\[
\phi^*\omega = \omega', {\;} \phi^* \weta = \weta' + du
\]
where $u \in \bO^0(U,\R^\times)$, defining the group $\cS(\omega)$ of $\lcs$ automorphisms of $(X,\omega)$. Let us write $\G_m$ for the subgroup $(1_X,u)$ of (conformal) rescalings. Brylinski's algebra $(\bO^*(X), \star)$ of differential forms (\cite{8}, see appendix {\bf a2}) has symplectic Hodge duality operator
\[
\alpha \wedge \star \beta := (\omega^{-\otimes k}(\alpha))(\beta) \cdot \frac{\omega^d}{d!}
\]
in degree $k = |\alpha| = |\beta|$, with $\omega^{-1} : T^\vee_X \cong T_X$ with $\vee$ for linear duals. \bigskip

{\bf 1.2} Working in a closed model category of differential commutative graded $\Q$-algebras, we have a commutative diagram ($\omega \mapsto \omega_2, {\;} \weta \to \omega_1$) 
\[
\xymatrix{
\OBlcs_\qt \ar@{.>}[d]^{t \to 0} \ar@{.>}[r] & \bL^*_{\Q[t]} \{\omega_1, \omega_2 {\; | \;} d\omega_2 = t \omega_1 \wedge \omega_2, {\;}, d\omega_1 = 0 \} \ar[d]^{- \otimes _\Q \Q[t=0]} \\
\OS_\Q \ar@{.>}[r] & \bL^*_\Q \{\omega_1, \omega_2 {\; | \;} d \omega_2 = 0, {\;} d \omega_1 = 0 \} {\;} ,}
\]
with dotted arrows signifying cofibrant (minimal model) replacements, following Neisendorfer \cite{54}(\S 5.6a) and Meier \cite{46}(\S 2.3).

I am indebted to Urs Schreiber for this -- these dcas are {\bf not} nilpotent \cite{36}[Def 1.10], but $\pi_1 = \Z$ has cohomological dimension one, which suffices for this construction -- and, more generally, for drawing my attention to the utility of rational-homotopy-theoretic methods, \eg \cite{23}(ex 2.41), \cite{27}(\S 4.6, open M5 branes), \cite{71} in mathematical physics. I also owe him thanks for his help and patience, but the mistakes here are my responsibility.\bigskip

{\bf 1.3.1} An $\lcs$ structure on $X$ defines a morphism $\omega_2 \to \omega_U, \omega_1 \to \weta_U$
\[
\omega_2^{-1}\OBlcs \to \bO^*_\R(X)
\]
of dcas, making the real differential forms on $X$ into an {\bf effective} $\OBlcs$ algebra.\bigskip

{\bf From here on $n = 2$}.

\newpage

{\bf1.3.2 Prop} {\it A collar $(0,1] \times Y \subset X$ for the boundary $\partial X = Y$ defines a commutative diagram 
\[
\xymatrix{
\OBlcs \ar[d] \ar[r] & \bO^*_\R(X) \ar[d] \\
\OS \otimes \Q[t] \ar[r] & \bO^*_\R(I \times Y) }
\]
presenting an $\bO^*\Omega S^2$ - algebra structure on $\bO^* Y$ as boundary data for an effective $\OBlcs$ structure on $\bO^*X$.}

[The right vertical map is defined by restriction. Note that $\bO^*$ signifies an algebra of differential forms, whereas $\Omega$ indicates a based loop space.] \bigskip

{\bf 1.3.3 } More precisely, fix an open book decomposition of $Y$, supporting a contact form ($\alpha \in \bO^1Y$ such that $\alpha \wedge d \alpha$ never vanishes, \cf {\bf a1}) \cite{18, 19, 20, 28}, \cite{29}(\S 11.1); then
\[
\omega_2 = -t^{-1}d(t \alpha), {\;} \omega_1 = t^{-1} dt = \weta  \in \bO^1 
\]
defines an effective $\bO^*B\lcs$ - algebra structure on $\bO^*(I \times Y)$. Indeed
\[
d\omega_2 = \weta \wedge \omega_2, {\; \rm and \:} \omega_2 \wedge \omega_2 = 2 \weta \wedge \alpha \wedge d\alpha \neq 0 {\;} .
\] 
The Kerr model \cite{11,48, 57}  is an interesting noncompact example.\bigskip

{\bf 1.4} See {\bf a1} for more on contact forms, but for the moment recall \cite{63} that a contact form has both an associated real $[\xi] \in \Pic_\R Y \cong H^1(Y,\Z_2)$ line bundle as well as the complex line $[\eta] \in H^2(Y,\Z) \cong \Pic_\C(Y)$ (related by $[\xi] \otimes_\R [\eta] = [\eta]$), and that an almost contact structure is a reduction of the structure group of $Y$ to 
\[
\T \rtimes \Z_2 \subset \SU(2) = {\rm Spin}(3) {\;} .
\]
The class $[\weta] \in H^1(Y,\R)$ can perhaps be regarded as a lift of $[\xi]$; {\bf NB} it is not the $\eta$ of \S 2.1 below. 

\bigskip

{\bf \S II Spin, contact, and {\bl}ine structures}\bigskip

{\bf 2.0} An almost contact structure on $Y$ as above is rationally indistinguishable from an $\bO^* \Omega S^2$ - algebra structure, but this loses two-torsion information. We argue here for working in a topos \cite{50}
\[
\{ X {\;} /{\;} |\Fin^\otimes_\C|^+ \}
\]
of spaces with enough basepoints -- alternatively, spaces endowed with a principal {\bl}ine or $\Omega^2 S^2 \simeq \Omega^2 \SU(2) \rtimes \Z$ bundle\begin{footnote}{a noncommutative harmonic oscillator}\end{footnote} -- as a rationally equivalent alternative. Such framings/orientations/generalized basepoints were proposed in \cite{49}(\S 2.3) as useful in understanding Cartier's renormalization flow in ${\rm MU}_\Q$-theory.

\newpage

{\bf 2.1} Recall that under Pontrjagin multiplication, the integral homology
\[
H_* (\Omega S^2,\Z) \cong \Z[\bx]
\]
is polynomial on a single generator $\bx$ of degree {\bf one}, and is therefore {\bf not} a (graded-) commutative algebra. The cohomology 
\[
H^*(\Omega S^2,\Z) \cong E([\omega_1]) \otimes P([\omega_2]) 
\]
{\bf is} however graded - commutative, with $[\omega_2]$ Kronecker dual to $\bx^2$. There is further background in appendix {\bf a1} below.

{\bf Definition} The map
\[
\C^\times \ni z \mapsto [q \to zqz^{-1}] \in \Maps_{\{1\}}(\bH^\times,\bH^\times)
\]
extends by identity maps when $z \in \{0,\infty\}$ to define 
\[
\bP^1(\C) = \C^\times \cup \{0,\infty\} \to \Maps_{{\rm deg} = 1}(S^3,S^3) {\;} ,
\]
(\ie if $z = r exp (i \alpha)$ with $r \in \R^\times$ then $z \to [q \to e^{i \alpha} q  e^{-i \alpha}]$, but if $r = 0,\infty$ it defines a Mercator-like compactification, with fixed points at the north and south poles).

Taking adjoints gives a loop map [{\bf a3}]
\[
[S^3 \to \Maps(\bP^1(\C),S^3) \cong \Omega^2 S^3] = \eta^2 \in \pi_5 S^3 \cong \Z_2 
\]
and thus a composition
\[
B \eta^2 : B\SU(2) \to \Omega S^3 \to \Omega S^2
\]
(the last map being the universal cover). The composition $\T \rtimes \Z_2 \to \SU(2) \to \Omega^2 S^2$ therefore assigns a $\bl$ine bundle to an almost contact structure via the remarkable relation
\[
\eta^3 = 4 \nu
\]
of Ravenel \cite{60}(1.1.14). Note that $\eta^2$ is the first nontrivial element in the cokernel of the stable $J$ - homomorphism.\bigskip

{\bf 2.2} A spin structure on $X$ thus has compatible associated $\bl$ine bundles
\[
\xymatrix{
Y \ar[r] \ar[d] & B\Spin(3) = B\SU(2) \ar[d] \ar[r]^-{B\eta^2} & B(\Omega^2 S^2) \ar[d] \\
X \ar[r] & B\Spin(4) =  B\SU(2)^{\times 2} \ar[r] & B(\Omega^2 S^2)^{\times 2} }
\]
on both $X$ and its boundary.

More generally, our topos 
\[
\{ X {\;} / {\;} |\Fin^\otimes_\C|^+ \} \to \{ X {\;} / {\;} B \Spin \}
\]
lies behind that of manifolds with spin structures, which already has a rich deformation theory \cite{38}. \newpage

{\bf 2.3 Proposition} {\it Compact Riemannian four-manifolds $(X,g,\omega)$ with effective $\bO^*B\lcs$ - algebra structures as morphisms, and their bounding three-manifolds (with compatible (\cf \S 1.2-3) $\bl$ine bundles) as objects, define a cobordism (topological, tri) category $\{Y \: /\: |\Fin^\otimes_\C|^+ \}_\lcs$.}

More precisely: its morphism objects
\[
\Mor_*(Y',Y) {\;\;} = \coprod_{X,\omega \supset \partial X = Y'_{\rm op} \sqcup Y} [\Met \dots /\cS(\omega)]
\]
are topological groupoids with 4D $\lcs$ cobordisms with $\bl$ine bundles on their boundaries as objects, and diffeomorphisms preserving these structures \begin{footnote}{\cite{24}(Ths 2,3),\cite{44}. Their geometric realizations or totalizations as simplicial spaces are countable \cite{43} CW spaces, not nilpotent in general}\end{footnote} as their morphisms. Such 4D cobordisms compose under glueing along compatible boundaries.

The effectiveness of an $\bO^*B\lcs$-structure permits the construction of a symplectic Hodge duality operator, and so \cite{2, 13, 41, 58, 65, 66} leads to an $\SL_2 \R$ - representation - valued Hodge-Lefschetz bidifferential cohomology theory \cite{2} with good finiteness properties, \cf {\bf a2} below. \bigskip

{\bf III An Introduction}\bigskip

{\bf \S 3.1} The points of a (? symplectic) manifold as a groupoid\\
{\bf \S 3.2} Euler, van Kampen, Mayer-Vietoris for groupoids\\
{\bf \S 3.3} compact subgroups of diffeomorphism groups\\
{\bf \S 3.4} Renormalization and Anderson localization\\
{\bf \S 3.5 Opening remarks} re GKM manifolds\bigskip

The narrative below replaces the previous $\S 3$; it attempts to explain why the classifying space (over $\Q$) for $\lcs$ structures might be interesting or useful. Behind that is the notion, for which I thank Urs Schreiber, that 
{\bf un}stable rational homotopy theory might be useful in the physics of diffeomorphism groups.

{\bf \S 3.1} the points of a manifold as a groupoid

{\bf 3.1.1} Suppose $X$ is a smooth compact connected manifold, either closed or with a connected boundary $\partial X = Y$; let  $\D(X)$ be its (\dots) Lie \cite{26, 56} group of diffeomorphisms (which take the germ of a collar of $Y$ to itself), and let $\D_0(X)$ be its identity component.

If we interpret $\D(X)$ as a group of permutations of $X$ then the transformation groupoid $[X // \D(X)]$ is the topological category $\cX$ of points of $X$ \cite{7, 68} \ie regarded as elements of the {\it smooth} manifold $X$. Similarly, given (for example) a symplectic form $\omega$ on $X$, there is a category $[X \sslash \cS(\omega)] := \cX_\cS$ of points of $X$ regarded as elements of a {\it symplectic} manifold (with group $\cS(\omega)$ of symplectic diffeomorphisms), and there is a forgetful functor
\[
[X \sslash \cS(X)] = \cX_\cS \to [X \sslash \D(X)] = \cX\;.
\]

Let us consider the  {\bf un}stable rational homotopy types defined by such objects, in terms of Sullivan - Quillen mixed (differential graded commutative) $\Q$ - algebras \cite{55, 62} which underlie equivariant Borel - de Rham cohomology \cite{3}, \ie Eilenberg-Mac Lane 
\[
H_\Q^*(X \times_\D E \D;\Q) := H_\D^*(X;\Q)
\]
cohomology of the simplicial totalization $|[X\ \sslash \D]| = X_{h\D}$ (which is an object of the weak homotopy type 
\cite{43} of) a CW space with countably many cells). 

We argue here that the geometric realization $X \mapsto |[X \sslash \D_0(X)]|$ of the category of points behaves much like a topological quantum field theory, with (the cohomology of) the space of points of a cobordism between three-manifolds being its bimodule transition operator in the sense of Dirac and Segal, \cf \S 3.5. The maximal compact toruses in such diffeomorphism groups are to some extent accessible \cite{12,39}, in particular when the cohomology is free over $H^*_T$ \cite{70}. \bigskip

{\bf 3.1.2} {\bf Motivation} for this comes from a topological category with {\bf locally conformally} symplectic ($\lcs$) manifolds $(X,\omega)$ as morphisms, their (\eg contact) boundaries being the objects of the category. Following work of Calabi, Banyaga, McDuff \cite{5, 9, 43}  and more recently Haller and Rybicki, the identity component $\cS_0(\omega)$ of the group
\[
\G_m \to \cS (X,\omega) \to \D(X) 
\]
of (lc)symplectic automorphisms $(u,\phi) \in \cS(X,\omega)$ of $X$ is an extension \cite{33}
\[
1 \to \cH \to \cS_0 \to \Phi \to 1
\]
of a perfect ($\cH = [\cH,\cH]$) group \cite{34}  of Hamiltonian (Lie $\cH \simeq C^\infty_0(X)$ \cite{8}(\S 3.4)) diffeomorphisms by a nilpotent Lie group $\Phi = \Phi(\omega)$ as in {\bf a5} below. \bigskip

{\bf 3.1.3} Let us write $\varepsilon : \Met(X) \to X$ for the space of (basepointed) Riemannian metrics $(g,x)$ on $X$. The action of $\D$ lifts to make $\varepsilon$ equivariant; it is a homotopy equivalence with cones for fibers, and we can consider, given $(X,\omega)$, the diagram of groupoids below (built mostly as a pullback of forgetful functors)
\[
\xymatrix{
[\lcs \times_{\G_m} \Met(X) \sslash \cH] \ar[d] \ar[r]^-\Phi & [ X \sslash \cS_0] \ar[d] \\
[\Met(X) \sslash \D_0]  \ar[r] & [ X \sslash \D_0]}
\]
to be the definition of its upper left corner, \ie a category of points of the space of Riemannian metrics on a manifold with a conformally compatible $\lcs$ structure. The lower left corner is a model (dim $X$ = 4) for the configuration space for Euclidean models of general relativity \cite{14, 42, 56}; the diagram attempts to model Euclidean gravity coupled to locally conformally symplectic matter of some \cite{11, 48, 57} sort. \bigskip

{\bf \S 3.2 Anderson localization} 

{\bf \S 3.2.1} Euler, Mayer-Vietoris, van Kampen \dots for groupoids

If instead we consider connected $X$ with {\bf two} boundary components, let 
\[
X :  Y'  \to Y_{-}, \; X' :  Y_{+} \to Y'' 
\]
denote a pair of manifolds $X,X'$, each with two boundary components $Y',Y_{-};Y_{+},Y''$, together with diffeomorphisms $\alpha_\pm : Y_\pm \cong Y$ to a reference manifold $Y$ carrying the germ of a normal bundle.
We can glue along $Y$ to define a coproduct manifold $X \circ X' : Y' \to Y''$ :
\[
\xymatrix{
Y \ar[d] \ar@{.>}[dr] \ar[r] & X \ar[d] \\
X' \ar[r] & X \circ_Y X' }
\]
again with two boundary components, as well as a smoothly embedded copy of $Y$ and a diffeomorphism $\alpha := \alpha_{-} \circ \alpha_+^{-1} : Y \to Y$ in $\D(Y) := \Delta$. If $\D(X \# X')$ denotes the group of diffeomorphisms of $X \circ X'$ compatible with $\alpha$ and the partition of the construction into left and right parts, then the forgetful functor
\[
[X \circ X' \sslash \D(X \# X')] \to [X \circ X' \sslash \D(X \circ X')] 
\]
defines an associative composition operation
\[
[X \sslash \D] \times_{[Y \sslash \Delta]} [X' \sslash \D'] \to [X \circ X' \sslash \D(X \circ X')] \;.
\]
We propose the mixed Quillen-Sullivan differential commutative algebra \cite{62}
\[
{\rm CE}^*(\D;\bO(X)) := \Lambda^*\cX
\]
(of Chevalley-Eilenberg cochains \cite{10}(Ch XIII \S 7) on $\D$ with values in $\bO^*(X)$) as candidates for a presentation
 \[
\Lambda^*\cY \to \Lambda^*(\cX \circ_\cY \cX')  \to \Lambda^*\cX \times \Lambda^*\cX'
\]
of these fiber products as HNN equalizers \cite{6}. In any case
\[
\{ X \to |[X \sslash \D(X)]| := \cX, \; X' \to |[X' \sslash\D(X')]| := \cX',\; Y \to |[Y \sslash \D(Y)]| := \cY \}
\]
\[
\mapsto  \cX \circ_\cY \cX' :=  [X \circ X' \sslash \D(X \circ X')] 
 \]
defines a functor
\[
(X,\partial X) \mapsto (\cX,\cY' \otimes \cY_{-} ) \mapsto  ( |\cX|_\Q,|\cY'|_\Q \times |\cY_{- }|_\Q) \;.
\] 
from a category with three-manifolds as objects and four-manifolds as morphisms, mapping in a motivated way to a category of homotopy types and their rational unstable dcas.

This deserves a nice presentation in 2-categorical terms which is outside my range. We'll work here instead with the homotopy groups of these objects, regarded as sheaves of modules over Spec $H^*_\D(X)$.]

{\bf 3.3.1} Recall that a topological (for example, connected) group $G$ acts by conjugation $g,K \to g \cdot K \cdot g^{-1}$ on the set
\[
\cpct_0 G := \{ K \subset G \;|\: K \; {\rm compact \; connected} \}
\]
of its connected compact subgroups, defining a functor 
\[
[\cpct_0 G \sslash G] \to [{\rm pt} \sslash G] \;.
\]
{\bf Definition} The object in the upper left corner 
\[
\xymatrix{
[X_\& \sslash \D_0] \ar@{.>}[d] \ar@{.>}[r] & [X \sslash \D_0] \ar[d] \\
[\cpct_0 \D \sslash \D] \ar[r] & [{\rm pt} \sslash \D] }
\]
defines Anderson localization $X_\& \to X$ (\ie away from the disordered part of a continuous medium \cite{70}).

{\bf 3.3.2} Following Quillen, $H_\D^*(X)$ defines a sheaf of graded algebras over Spec $H^{2*}_\D({\rm pt})$, but the $\D$-equivariant sheaves defined by restriction to subgroups $K$ of $\D$ will be of interest, and we will write 
\[
H^*_{\& \D}(X) :=  [K \mapsto H_K^*(X)]  
\]
for the sheaf of abelian groups over $[\cpct_0 \D_0 \sslash \D_0]$ so defined \begin{footnote}{rather than pass to some limit we keep the hood open \cite{67} while working}\end{footnote}.

{\bf 3.3.3} In the theory of finite groups, the Atiyah-Quillen stratification \cite{59} of Spec $H^*_\Gamma ({\rm pt} ; \F_p)$ is defined by the directed category $\AQ_p(\Gamma)$ with elementary abelian $p$-subgroups 
\[
\F_p^k \cong A \subset \Gamma
\] 
of $\Gamma$ as objects, and morphisms
\[
\Mor(A,A') := \{ \phi \in \hom(A,A') \; | \; \exists \gamma \in \Gamma, \; \phi(a) = \gamma \cdot a \cdot \gamma^{-1} \} \;.
\]
Quillen shows that the induced homomorphism 
\[
H^*_\Gamma ({\rm pt}) ; \F_p) \to {\rm colim}_{A \in \AQ} H^*_A ({\rm pt}) ; \F_p) := H^*_{\AQ \Gamma} ({\rm pt};\F_p)
\]
represents a topological (\ie modulo nilpotents) isomorphism of the underlying schemes.

In our case, the natural transformation
 \[
H^*_\D \to {\rm colim}_{K \in \cpct_0 D} H^*_K := H^*_{\& \D}
\] 
suggests restriction to a subspace of Spec $H^{2*}_\D$. For reasons of time we restrict to similar colimits over compact connected abelian {\bf torus}es.\bigskip 

{\bf 3.4} The isotropy group $\Stab(g)$ of a compact Riemannian manifold is compact \cite{53}, so its identity component defines a functor 
\[
 g \mapsto \Stab(g)_0 : [\Met(X) \sslash \D(X)] \to [\cpct_0 \D \sslash \D] \;.
\]

{\bf Proposition} {\it The fiber product
\[
\xymatrix{
[\Met (X) \sslash \D_0] \ar@{.>}[d] \ar@{.>}[r] & [X \sslash \D_0] \ar[d] \\
[\cpct_0 \D_0 \sslash \D_0] \ar[r] & [{\rm pt} \sslash \D_0] }
\]
defines a consequent isomorphism}
\[
H^*_{\D_0}(\Met (X)) \cong H^*_{\D_0}(X_\&)
\]

{\bf `proof':} a hypothetical Leray sseq 
\[
H^*(\cpct_0 \D_0/ \D_0, {\mathcal H}^*_\Stab) \Rightarrow H^*_{\D_0} (\Met(X))
\]
collapses $\Box$ .

In this language we rewrite the diagram of \S1.2 as 
\[
\xymatrix{
[X_\& \slash \cH] \ar[d] \ar[r]^-{\Phi_\&} & [X \sslash \cS_0] \ar[d] \\
[X_\& \sslash \D_0]  \ar[r] & [X \sslash \D_0]}
\]
with the left side considered as a cofibrant replacement or regularization of the right side, interpreting gravity as renormalization \cite{49}. \bigskip

{\bf \S 3.5 Opening remarks}\bigskip

The characteristic zero cohomological behavior of compact Lie groups is well-understood in terms of their maximal tori \cite{31}(\S 8.4), \cite{39}. Let us write $\cA(G)$ for the collection of compact connected abelian subgroups of $G$, and let
us regard the functor 
\[
\cH_\cA (X) : [\cA(\D_0) \sslash \D_0]  \ni T \mapsto  H^*_T(X) \in H^*_{\D_0} - (\Mod)
\]
as shorthand for the cohomological study of torus actions on manifolds; for simplicity let us restrict to equivariantly formal (\ie free over $H^*_\T({\rm pt})$) actions \cite{31}(Th 1.2.2).

In fact each vertex of the diagram above has its own relevant category and literature of torus actions, and the post to which this is an introduction was motivated in particular by the Hodge-Lefschetz symplectic cohomology \cite{2, 8, 13, 41, 58, 65, 66} associated to the upper right entry. The lower left adds characteristic class technology \cite{38} to the toolbox. A great deal is known about Hamiltonian circle actions on symplectic four-manifolds \cite{1,4,47},\dots, though alas not by me. Let us leave this here following the account \cite{32} of Guillemin and Zara:

{\bf Example} [GKM \S 7] a $T$-space $X = \bigcup_{i \geq 0} \sk_i X$  is the union of its ($T$ -space) skeleta (ordered by stabilizer codimension), so its one-skeleton consists of fixed points and circles. If these are finitely many and the cohomology is free over $H^*_\T$ then $\sk_1X/T$ is a (decorated \cite{30}) connected graph which determines Spec $H^{2*}_T(X)$ as a topological space.
 
This is the basis for the analogy between such GKM manifolds and (physical issues surrounding the language of)
  Feynman graphs. The one-skeleton $\sk_1 X$ is built of dipole bubble collapse maps 
\[
[\C P^1 \sslash \T]  \to (\C P^1/\T) (\simeq [(\C P^1/\T) \sslash 1]) \cong [-1,+1])
\]
defining the half-edges $[\C_{\leq 1} \sslash \T] :\simeq c^{-1}$ of its graph as in \cite{51}(\S 2.4). 

\newpage

{\bf appendices}\bigskip

{\bf a1} In these notes $Y \cong Y(\bB)$ is in principle a compact connected closed oriented three-manifold, presented as an open book 
\[
\xymatrix{
Y \supset Y - \bB \ar[r]^\Sigma & S^1 }
\]
with a punctured surface $\Sigma$ as leaf, bound by a link $\bB \cong {\rm B} \times S^1$, \cite{18}, \eg the mapping torus of a two-manifold. Any (\dots) three-manifold has such presentations. 

[Thurston-Winkelnkemper] \cite{64} There is a contact form $\alpha, \ie {\;} \alpha \wedge d \alpha \neq 0$ supported by $\bB$. 

[Giroux] \cite{28} : The space of contact forms supported by a binding $\bB$ is connected (under isotopy = smooth deformation of the contact form $\alpha$\begin{footnote}{accomplished by stabilizing \cite{18}(\S 4.31 Fig 12] under certain (? 2-categorical) moves}\end{footnote}).

Taking universal covers on the right leads to the diagram
\[
\xymatrix
{Y \ar[d]^\cra & \ar[l] Y - \bB \cong \Sigma \times_\Z \R \ar[d] & \ar[l] \Sigma \times \R \ar[d] \\
\Omega S^2 & \ar[l] \Omega S^2 & \ar[l]  \Omega \SU(2) }
\]
and thus by compactification to $a : \Sigma_+ \to \Omega^2 S^3$ inducing
\[
H_2(\Sigma_+,\Z) \to H_2(\Omega^2 S^3,\Z) = \Z_2
\]
[recall $H_{*>1} = \oplus \beta_p H_{1+*} (\Omega^2 S^3,\Z_p)$ (with the latter polynomial on generators of degree $2(p^k - 1)$ or so), defining a class in $H^2(\Sigma_+,\Z_2)$ perhaps related to the Hurewicz image of $\eta^2$?]

Poincar\'e's lemma in link calculus \cite{29}(Prop 4.5.11, Exercises 5.3.13] associates to (the link $\bB$ binding) $Y(\bB)$, a simply-connected four-manifold $X \supset \partial X \cong Y$ with intersection form $H^2(X/Y;\Z)$ defined on $H_2(X,\Z)$ as in \cite{28} (\S 5.4). The coboundary
\[
H^1 Y \to H^2 X/Y : \Pic_\wzw Y \to \Pic_\C X/Y
\]
sends $\bl$ines on the boundary to complex lines on the interior of $X$.\bigskip

{\bf a2} Following \cite{13}(\S 2), \cite{61}(\S 7), the real Lie algebra $\SL_2(\R)$ acts on the ungraded vector space $\bO^\bullet(X)$ defined by the hyperplane intersection $ L = \omega \wedge - $ operator and its symplectic adjoint $L^\star = - \star L \star$, with 
\[
[L,L^\star] = \Sigma_{k \geq 0} {\;}  (2n - k) {\;}  \Pi_k
\]
$\Pi_k$ being projection to grade $k$. With the symplectic coboundary operator
\[
\delta_k = (-1)^k \star d_k \star  = d_{k-1} L^\star - L^\star d_k
\]
then [AOT lemma 1.3]  for $k \in \N$ we have operators $\fd_k = d_k - k \weta$ with 
\[
\fd_k^2 = 0, {\;} \delta_k \delta_{k+1} = 0, {\;} \delta_k \fd_k + \fd_{k-1} \delta_k = 0 
\]
defining a {\bf bi}differential $\Z$-graded Hodge-Lefschetz dca $\bO^*_\HL(X)$ with well-behaved associated elliptic complexes when $X$ is compact \cite{2}(Prop 2.2), \cite{65,66}.

{\bf remark} The circle $\T$ acts on the right on ${\rm Sl}_2(\R)$, and thus on the product 
\[
T_{{\rm Sl}_2(\R)} \times_{\SL_2(\R)} \bO^\bullet(X) \to {\rm Sl}_2(\R)
 \]
of its tangent bundle with $\bO^\bullet(X)$ to define a bundle of graded dcas over the quotient disk ${\rm Sl}_2(\R)/\T$, which is reminiscent of the cohomology 
\[
\Spec H^{2*}(V) \to \Spec H^{2*}(\C P^\infty) \cong Spec {\;} \Q[\varkappa] = \A^1_\Z
\]
of a projective variety $V \to \C P^N$, perhaps regarded as the upper half-plane.\bigskip
 
{\bf a3} Verification of \S 2.1: The diagram
\[
\xymatrix{
S^3 \times S^3 \ar[r]^-{\eta^2 \times \eta^2} \ar[d]^{\mu_{S^3}} & \Omega^2 S^3 \times \Omega^2 S^3 \ar[d]^{\mu_{\Omega^2S^3}} \\
S^3 \ar[r]^{\eta^2} & \Omega^2 S^3}
\]
commutes, because across the top we have
\[
q,k \to [z \to zqz^{-1}],[w \to wkw^{-1}] {\;},
\]
then down along the right we have 
\[
 [z \to zqz^{-1}],[w \to wkw^{-1}] \mapsto [u \to u \times u \to uqu^{-1} \cdot uku^{-1} = ukqu^{-1}
 \]
 while around the left side we have $q,k \mapsto [u \to uqku^{-1}]$ {\;} .
 
 To be more precise about basepoints:  in somewhat mixed notation, we have a commutative diagram
 \[
 \xymatrix{
 \{0,\infty\} \times (\T \times \R^3_+) \ar[d]^\subset \ar[r]^-{{\rm pr}_1} &  \{0,\infty\} \ar[d]^\subset \\
 [0,\infty] \times (\T \times \R^3_+) \ar[d] \ar[r]^-\rho & \R^3_+ \ar[d]^\cong \\
 (S^2 = \bP^1\C) \times S^3 \ar[r]^-\brh & S^3 }
 \]
 in which $\rho(r,\omega,q) = \omega q \omega^{-1}$ for $r \in \R^\times$ and $\rho(0,-,-) = 0, \rho(\infty,-,-) = \infty$. The downward left collapse $\{0,\infty\} \times \T \to [0,\infty] \times \T \to S^2$ is a cofibration but the right side is just an inclusion. The assertion is not that $\rho$ is continuous, but that collapsing the two endpoints of a cylinder does define a continuous if not smooth $\brh$.\bigskip

{\bf a4} This is a remark about a rough circle of ideas, a homotopy-theoretic version of Wess-Zumino-Witten theory, where 

{\it A principal $\Omega^2 S^2$ (or $\bl$ine) bundle on $X$ is homotopically equivalent to a field $X \to |{\rm Fin}_\C^\otimes |^+$ of (virtual) approximations of the cx line by finite subsets.}

 The set of such equivalence classes has group structure, defining the homomorphism
\[
\pi^2(S^1 \wedge X) := \Pic_\wzw (X) \to \Pic_\C (X) \cong H^2(X,\Z) 
\]
induced by a generator of $H^2(\Omega S^2,\Z)$.

A $\bl$ine bundle has an underlying (cx) line bundle; they are candidates for homotopy-theoretic gauge fields, in particular on 3-manifolds. They can be interpreted as fields of WZW models, as follows:

Taking loopspaces sends fibrations to fibrations; applied to Hopf
\[
\T = S^1 \to \SU(2) = S^3 \to \C P^1 = S^2 
\]
we get the universal cover 
\[
\Omega \T \simeq \Z \to \Omega \SU(2) \to \Omega S^2
\]
so looping once again defines a homotopy equivalence 
\[
\Omega^2 \SU(2) {\; } (\cong \Omega^2 \R^3_+) \to \Omega^2_e S^2
\]
with the identity component of a semi-direct product
\[
\Omega^2 S^2 \simeq  \Omega^2 \SU(2) \rtimes \Z \sim [\Omega^2 \SU(2) \sslash \Z] 
\]
though one must be careful because group-like spaces are not necessarily groups.

Note that  $\Z \cong \pi_0 (\Omega^2 S^2)$ builds the writhe or determinant of a braid (as in \S  8 of \cite{52}:  a map $X \to Y$ of Poincar\'e duality spaces, with a lift  of its normal spherical bundle to the classifying space\begin{footnote}{which is the (to me mysterious) universal cover $\Omega^2 \SU(2) \langle$ }\end{footnote} of (degree one self-maps/homotopy automorphisms of) $S^3$) is Steenrod representable, by work of Hopkins and Mahowald) into geometry. In this language there's a topos of (Poincar\'e duality) spaces over $B\HAut(S^3)$ which is arguably where solutions of variational geometric problems live, while this note argues/claims that the topos of (4D) spaces over $B\HAut(S^2)$ is a good place for renormalizeable 4D Riemannian geometry/qEuclidean GR. 

There is a nice 3-form 
\[
\ga d \ga + \lambda \cdot \kappa(\ga,\ga,\ga)
\]
on the space of $\mathfrak{su}(2)$ - connections, where  $\kappa$ is related to the Killing form of Lie theory and to the generator of $\pi_3$ of a (...) simple Lie group; it leads to invariants of three-manifolds defined without reference to Riemannian geometry. The Bott-Cartan element underlying $\kappa$ is the class of origin for $4 \nu_{(2)} = \eta^3 \in \pi^S_3$ in the stable homotopy groups of spheres; note that both $\eta$ and $\eta^3$ are in the image of the $J$-homomorphism, but $\eta^2$ is not. \bigskip

{\bf a5} notes on Haller and Rybicki  \cite{33}(\S (4.3, 5.6, 6.3):

$\cS$ is filtered as follows :
\[
\Phi :  \cS_0/\cS_1 \cong \Delta, {\;} \cS_0 \to H^0_c/\Delta
\]
\[
\Psi :  \cS_1/\cS_2  \cong  \Gamma, {\;} \cS_2 = {\rm ker} \Phi  \to H^1_d /\Gamma
\]
\[
R :  \cS_2/\cS_3  \cong \Lambda, {\;} \cS_ = {\rm ker} \Phi  \to H^{2n}_{d^{n+1}}/\Lambda
\]
with $\cS_3 = {\rm ker \;} R$. \bigskip

{\bf a6 extract} from a letter to J Baez : 

\dots thinking again about the Kerr model. It's easy to define and there are interesting things that are easy to say about it. It's a fundamental example in GR and could/should/ought to be understood better in topology: \bigskip

{\bf 1} You can write the Kerr-Schild coordinate system on the back of an envelope:
\[
{\sf ks} : \R^2 \times S^2 \ni (r,z,\theta,\phi) \mapsto (t,x,y,z) \in \R \times (\R^3_+ - \T_a) 
\]
where 
\[
 \T_a = \{(x,y,0) {\; | \;} x^2 + y^2 = a^2 \} \subset \R \times \R^3_+
\]
is the circle of radius $a$, from one end of time to the other \begin{footnote}{Usually the coordinates are $u_\pm = r \pm t \in \R^2 (\times S^2$); we restrict to the diagonal $u_+ = u_- = r$, ie $t = 0$ and interpret the spatial slice as constant in time.}\end{footnote}, and
\[
x + iy : = (r- ia) \sin \theta \cdot  e^{i\phi}  {\;} {\rm and \;} z : = r \cos \theta {\;};
\]
note $r \in \R$, and that if $\lambda := du_+  + \Theta d\phi$ with $\Theta = a t \sin^2 \theta$ then 
\[
d(t\lambda) = 1/2 du_+ \wedge du_- + d \Theta \wedge d\phi 
\]
seems to be a symplectic two-form \cite{48}(\S 2.3).\bigskip

{\bf 2} On this coordinate patch the Kerr (Ricci-flat Lorentz) metric (cf wikipedia), a very plausible model for a rotating black hole, is defined. The following is concerned entirely with the elementary topology of this manifold and some of its compactifications, but not at all with its metric. 

Topologically the domain and range of $\sf ks$ are unlike: the target is a $K(\pi,1)$ space for $\pi = \Z$ (indeed an $(\R \times (\R^2_+ - \infty)$)-bundle over the circle), while the domain is simply-connected. Note the interesting $\{\pm\}^5/\{\pm\} \cong \Z_2^4$ - action. This deserves graphics, \eg the circle of radius $a$ as the unknot in $\R^3$ with a pinched point on its normal bundle as below. \bigskip

{\bf 3} We observe that
\[
x^2 + y^2 + z^2 = a^2 (1 - (r^{-1}z)^2) {\;} .
\]
Let us simplify for the moment by taking $y=0$, so that
\[
(x^2  +  z^2) = a^2 (1 - (r^{-1}z)^2)
\]
defines the quartic surface 
\[
\kerr := \{[r:x:z:a] {\; | \:} r^2(x^2 + z^2) = a^2 (r^2 - z^2) \} \subset \bP^3(\R)
\]
in real projective three-space.\bigskip

{\bf 4 Proposition} {\it There is a classical Lefschetz pencil
 \[
 \xymatrix{
\kerr \ar[d]^{- \log |\tan \theta/2|} \ar[r]^-{\sf ks} & \bP^3(\R) \ni [r:x:z:a ] \ar[d]^\zeta \\
\bP^1(\R) \ar[r]^-{\rm arctanh} & \T/\pm \ni [z: r]} 
\]
of one-dimensional compact real algebraic curves, which fibers $\kerr$ by ellipses} :

of eccentricity $e = \sin \theta$ in the $(\varkappa,\rho)$)-plane, presenting $\kerr$ as a torus $\T^2/\T \times 0$ with a pinched cycle; alternatively, as a trivial circle bundle with a degenerate fiber, or as a surface $\R^2_+/\{0,\infty\}$ of genus zero with a degenerate double point \dots \bigskip

{\bf 5} There are descriptions from many viewpoints.  This compactification of the Kerr model looks topologically like the usual hyperbolic cone but differs metrically; it's steeper, quartic rather than quadratic. Relaxing the $y=0$ condition blows up the points of $\kerr$ to circles of revolution in the domain of $\sf ks$, which are the celestial spheres above the points of the $(\varkappa,\rho)$ - plane. 
\bigskip

{\bf 6 Proof} : First of all, $\kerr |_{a = 0}$  is the (very degenerate) split quartic \cite{9}
\[
r^2 (x^2 + z^2) = 0
\]
This calls for a genuine algebraic geometer: the classical case $xy = uv$ is classical \cite{9} but I don't understand the squares, {\bf tbc}\bigskip

Otherwise, let
\[
\varkappa : = a^{-1} x, {\;} \rho : = a^{-1} r, {\;} \zeta : = a^{-1}z {\;}.
\]
Then we have
\[
\rho^2 (\varkappa^2 + \zeta^2) = \rho^2 - \zeta^2
\]
so
\[
\rho^2 (\varkappa^2 - 1)  + (1 + \rho^2) \zeta^2 = 0 {\;},
\]
hence
\[
\zeta^2 = (1 + \rho^{-2})^{-1} (1 - \varkappa^2) {\;}.
\] 

{\bf 7} This pencil follows the azimuthal coordinate $\zeta$ but funny things happen at $\zeta = 0$ and it is useful to introduce a parameter $\tz = \cos \theta = \rho^{-1}\zeta$: then the fiber 
\[
(1 + \rho^2)\tz^2 = 1- \varkappa^2 {\;} \Rightarrow {\:} \varkappa^2 + \tz^2 \rho^2 = 1 -  \tz^2
\]

above $[z:r] = \zeta$ defines an ellipse of eccentricity 
\[
e^2 = 1 - \tz^2
\]
in the $(\varkappa,\rho)$-plane. \bigskip

{\bf 8} Composition around the left side of the diagram is a classical trigonometry exercise in showing that $\zeta = \rho \cos \theta$. This defines a nice map to $\bP^1(\R)$ with which to construct a universal simply-connected cover for $\kerr$. Homotopically, it's a string of spheres - balloons in Milnor, a bubble-chain of $\C P^1$s attached at $\{0,\infty\}$ : here covered by strings of 3-spheres blown up from $\kerr$ by the SO(2) -symmetry in $(x,y)$. It's not the classical chain of Friedman \ae ons. \dots 

\bibliographystyle{amsplain}

\end{document}